\documentclass[preprint,12pt]{elsarticle}
\usepackage{amsmath}
\usepackage{amssymb}

\input{tcilatex}

\begin{document}

\begin{frontmatter}



\title{
On the maximal reduction of games
}

\author{Monica Patriche}

\address{
University of Bucharest
Faculty of Mathematics and Computer Science
    
14 Academiei Street
   
 010014 Bucharest, 
Romania
    
monica.patriche@yahoo.com }

\begin{abstract}
We study the conditions under which the iterated elimination  of strictly dominated strategies 
is order independent and we identify a class of discontinuous games for which order does not
 matter. In this way, we answer the open problem raised by M. Dufwenberg and M. Stegeman 
(2002) and generalize their main results. We also establish new theorems concerning the 
existence and uniqueness of the maximal game reduction when the pure strategies are dominated 
by mixed strategies. 

\end{abstract}

\begin{keyword}
Game theory, \
strict dominance,\
iterated elimination, \
order independence, \
maximal reduction.\


\end{keyword}

\end{frontmatter}



\label{}





\bibliographystyle{elsarticle-num}
\bibliography{<your-bib-database>}







\section{Introduction}

The question raised by Pearce (1984), concerning the rationalizable
strategic behaviour of the players in noncooperative strategic situations
was followed by a great amount of literature. It seemed to attract the
interest of researchers from Game theory. The first step of research in this
area was made by Bernheim (1984) and Pearce (1984), who defined the
rationalizable strategies of a strategic game by using iterative processes
of elimination of dominated strategies that were considered 'undesirable'.

This procedure led to the issue of order independence, which was studied by
many authors. They searched for classes of games and defined dominance
relations under which the result of the iterative process of removal of
dominated strategies does not depend on the order of removal.

Gilboa, Kalai and Zemel (1990) provided conditions (including strict
dominance) which guarantee the uniqueness of the reduced games. Marx and
Swinkels (1997) defined nice weak dominance, proved that under this order
relation, order does not matter. The main result of Dufwenberg and Stegeman
(2002) concerns a class of games for which a unique and nonempty maximal
reduction exists. The properties satisfied by games for which the iterated
elimination for strictly dominated strategies (IESDS) preserves the set of
Nash equilibria are the compactness of the strategy spaces and the
continuity of payoff functions. The authors also proved that if, in
addition, the payoff functions are upper semicontinuous in own strategies,
then the order does not matter. Chen, Long and Luo (2007) provided a new
definition of IESDS that proved to be suitable for all types of games and
also order-independent. Apt's approach (2007) uses operators on complete
latice and their transfinite iterations. The monotonicity of the operators
assures the order independence of iterated eliminations. Apt's paper (2007)
provides an analysis of different ways of iterated eliminations of
strategies. The notions of dominance and rationalizability are involved by
other two strategy elimination procedures studied by Apt (2005). In order to
study the problem of order independence for rationalizability, the author
considers three reduction relations on games and belief structures.

In this paper, we identify a class of discontinuous games for which order
independence holds, generalizing the main results of Dufwenberg and Stegeman
(2002). The payoff functions are transfer weakly upper continuous in the
sense of Tian and JZhou (1995). These authors defined the transfer upper
continuity and proved generalizations of Weierstrass and of the maximum
theorem. We also establish results for game reductions in which the pure
strategies are dominated by mixed strategies. We use some notions of
measurability and especially some results of Robson (1990).

The paper is organized in the following way: Section 2 contains
preliminaries and notations. Generalizations of Dufwenberg-Stegeman Lemma
are presented in Section 3. The mixed strategy case is treated in Section 4.
The concluding remarks follow at the end.

\section{Preliminaries}

Dufwenberg and Stegeman (2002) concluded that it remained an open problem to
identify classes of games for which order independence holds, outside of the
compact and continuous class.

We are searching to solve this problem. In order to reach this aim, we first
introduce the notions of games, parings, dominance and game reduction,
following that, in the next subsection, we discuss the transfer upper
continuity, a concept due to Tian and Zhou, which characterizes the payoff
functions of a class of games which generalizes than that one of Dufwenberg
and Stegeman. In section 3, we will prove that, in this case the iterated
elimination of strictly dominated strategies (IESDS) also produces a unique
maximal reduction.

\subsection{Games, Parings, Dominance and Reduction}

In the paper called "Equilibrium points in n-person games" (1950), Nash
describes without formalizing, the concepts of the n-person game and the
equilibrium of the attached game. He defines the n-person game, where each
player has a finite number of strategies and each n-tuple of strategies
corresponding to a given set of players wins. Any n-tuple of strategies can
be regarded as a point in the product space of sets of players' strategies.
A point of equilibrium is an n-tuple of strategies such that every player's
strategy brings the maximum payout for that player, against n-1 strategies
of the other ones.

We give the formal definition of an n-person game below.\medskip\ 

\begin{definition}
The normal form of an $n$-person game is $G=(I,(G_{i})_{i\in
I},(r_{i})_{i\in I})$, where, for each $i\in I=\{1,2,...,n\}$, $G_{i}$ is a
non-empty set (the set of individual strategies of player $i$) and $r_{i}$
is the preference relation on $\prod\nolimits_{i\in I}G_{i}$ of player $i$.
\end{definition}

The individual preferences $r_{i}$ are often represented by utility
functions, i.e. for each $i\in \{1,2,...,n\}$ there exists a real valued
function $u_{i}:\prod\nolimits_{i\in I}G_{i}\rightarrow \mathbb{R}$ (called
the utility function of $i$), such that $xr_{i}y\Leftrightarrow u_{i}(x)\geq
u_{i}(y),\forall x,y\in \prod\nolimits_{i\in I}G_{i}.$

Then the normal form of n-person game is $(I,(G_{i})_{i\in I},(u_{i})_{i\in
I})$.

\textit{Notation. }Denote $x_{-i}$=$(x_{1},...,x_{i-1},x_{i+1},...,x_{n})$
and $G_{-i}=\prod\nolimits_{j\in I\backslash \{i\}}G_{j}.$

\begin{definition}
The\textit{\ Nash} \textit{equilibrium }for the game $(I,(G_{i})_{i\in
I},(u_{i})_{i\in I})$ is a point $x^{\ast }\in \prod\nolimits_{i\in I}G_{i}$
which satisfies for each $i\in \{1,2,...,n\}:u_{i}(x^{\ast })\geq
u_{i}(x_{-i}^{\ast },x_{i})$ for each $x_{i}\in G_{i}.$
\end{definition}

Further we will assume that for each $i\in I,$ the set $G_{i}$ is a
Hausdorff topological space and $\prod\nolimits_{i\in I}G_{i}$ is endowed
with the product topology.

\begin{definition}
The game $G$ is called
\end{definition}

i) compact if $G_{i}$ is compact for each $i\in I;$

ii) own-uppersemicontinuous if $u_{i}(\cdot ,s_{-i})$ is upper
semicontinuous for each $i\in I$ and for each $s_{-i}\in G_{-i};$

iii) continuous if $u_{i}$ is continuous for each $i\in I.$

\begin{definition}
(Dufwenberg and Stegeman, 2002). A paring of $G$ is a triple $%
H=(I,(H_{i})_{i\in I},(u_{i}^{\prime })_{i\in I}),$ where $H_{i}\subseteq
G_{i}$ and $u_{i}^{\prime }=u_{i|\prod\nolimits_{i\in I}H_{i}}.$
\end{definition}

A pairing is nonempty if $H_{i}\neq \emptyset $ for each $i\in I.$

\begin{definition}
Given a pairing $H$ of $G$, the strict dominance relation $\succ _{H}$ on $%
G_{i}$ can be defined$:$
\end{definition}

for $x,y\in G_{i},$ $y\succ _{H}x$ if $H_{-i}\neq \emptyset $ and $%
u_{i}(y,s_{-i})>u_{i}(x,s_{-i})$ for each $s_{-i}\in H_{-i}.$

\begin{remark}
$\succ _{H}$is transitive.
\end{remark}

Let us consider parings $G,H$ with the property that $H_{i}\subseteq G_{i}$
for each $i\in I.$ We give here the definition of game reduction used by
Dufwenberg and Stegeman (2002), in order to generalize their main results,
following that, in the next section we will introduce other types of
reduction and discuss the relationships amongst them.

\begin{definition}
i) $G\rightarrow H$ is called a reduction if for each $x\in G_{i}\backslash
H_{i}$, there exists $y\in G_{i}$ such that $y\succ _{G}x.$
\end{definition}

ii) the reduction $G\rightarrow H$ is called fast if $y\succ _{G}x$ for some 
$x,y\in G_{i}$ implies $x\notin H_{i}.$

iii) the reduction $G\rightarrow ^{\ast }H$ is defined by the existence of
(finite or countable infinite) sequence of parings $A^{t}$ of $G,$ $%
t=0,1,2...$, such that $A^{0}=G,$ $A^{t}\rightarrow A^{t+1}$ for each $t\geq
0$ and $H_{i}=\cap _{t}A_{i}^{t}$ for each $i\in I;$

iv) $H$ is said to be a maximal $(\rightarrow ^{\ast })$-reduction of $G$ if 
$G\rightarrow ^{\ast }H$ and $H\rightarrow H^{\prime }$ only for $%
H=H^{\prime }.$

\subsection{Transfer upper continuity}

Tian and Zhou \ (1995)relaxed the continuity assumptions on functions and
correspondences which can be used in some economic models. Their work was
motivated by questions concerning the minimal conditions under which a
function reaches its maximum on a compact set or the set of maximum points
of a function defined on a compact set is non-empty and compact. Tian and
Zhou introduced the transfer continuities and generalized the Weierstrass
Theorem by giving a necessary and sufficient condition for a function f to
reach its maximum on a compact set.

We are providing here the concepts of transfer upper semicontinuity and
transfer weakly upper continuity for functions, the concept of transfer
closed-valuedness for correspondences and some of their properties.\medskip

Let $X,Y$ be subsets of topological spaces.

\begin{definition}
A function $f:X\rightarrow \mathbb{R}$ is said to be upper semicontinuous on 
$X$ if $\{x\in X:f(x)\geq r\}$ is closed in $X$ for all $r\in \mathbb{R}.$
\end{definition}

\begin{definition}
(Tian and Zhou, 1995) A function $f:X\rightarrow \mathbb{R\cup \{-\infty \}}$
is said to be transfer upper continuous on $X$ if for points $x,y\in X,$ $%
f(y)<f(x)$ implies that there exists a point $x^{\prime }\in X$ and a
neighborhood $\mathcal{N}(y)$ of $y$ such that $f(z)<f(x^{\prime })$ for all 
$z\in \mathcal{N}(y).$
\end{definition}

\begin{definition}
(Tian and Zhou, 1995) A correspondence $F:X\rightarrow 2^{Y}$ is said to be
transfer closed-valued on $X$ if for every $x\in X,$ $y\notin F(x)$ implies
that there exists $x^{\prime }\in X$ such that $y\notin $cl$F(x^{\prime }).$
\end{definition}

\begin{remark}
(Tian and Zhou, 1995) It is clear that, for any function $f:X\rightarrow 
\mathbb{R\cup \{-\infty \}}$, the correspondence $F:X\rightarrow 2^{X}$
defined by $F(x)=\{y\in X:f(y)\geq f(x)\}$ for all $x\in X$ is transfer
closed-valued on $X$ if and only if $f$ is transfer upper continuous on $X.$
\end{remark}

The next lemma characterizes the correspondences which have transfer
closed-values.

\begin{lemma}
(Tian and Zhou, 1995) Let $X$ and $Y$ be two topological spaces, and let $%
F:X\rightarrow 2^{Y}$ be a correspondence. Then, $\cap _{x\in X}$cl$%
F(x)=\cap _{x\in X}F(x)$ if and only if $F$ is transfer closed-valued on $X.$
\end{lemma}

The next property is a necessary condition for a function to have a maximum
on a choice set $G.$

\begin{definition}
(Tian and Zhou, 1995) A function $f:X\rightarrow \mathbb{R\cup \{-\infty \}}$
is said to be transfer weakly upper continuous on $X$ if, for points $x,y\in
X,$ $f(y)<f(x)$ implies that there exists a point $x^{\prime }\in X$ and a
neighbourhood $\mathcal{N}(y)$ of $y,$ such that $f(z)\leq f(x^{\prime })$
for all $z\in \mathcal{N}(y).$
\end{definition}

Theorem 1 generalizes the Weierstrass theorem.

\begin{theorem}
(Tian and Zhou, 1995) Let $X$ be a compact subset of a topological space and
let $f:X\rightarrow \mathbb{R\cup \{-\infty \}}$ be a function. Then $f$
reaches its maximum on $X$ if and only if $f$ is transfer weakly upper
continuous on $X.$
\end{theorem}

Morgan and Scalzo (2007) defined the upper pseudocontinuity and proved the
existence of Nash equilibrium for economic models with payoff functions
having this property.

\begin{definition}
(Morgan and Scalzo, 2007) Let $X$ be a topological space and $f:X\rightarrow 
\mathbb{R}$. $f$ is said to be upper pseudocontinuous at $z_{0}\in X$ such
that $f(z_{0})<f(z),$ we have $\lim \sup_{y\rightarrow z_{0}}f(y)<f(z_{0}).$
\end{definition}

\begin{remark}
The class of upper pseudocontinuous functions is strictly included in the
class of transfer upper continuous functions introduced by Tian and Zhou.
\end{remark}

\section{Generalizations of Dufwenberg-Stegeman Lemma}

The following lemma is due to Dufwenberg and Stegeman (2002).

\begin{lemma}
If $G\rightarrow ^{\ast }H$ for some compact and own-uppersemicontinuous
game $G,$ and $y\succ _{H}x$ for some $x,y\in G_{i}$ and $i\in I,$ then
there exists $z^{\ast }\in H_{i}$ such that $z\nsucc _{H}z^{\ast }\succ
_{H}x $ for each $z\in G_{i}.$
\end{lemma}

Let $G\rightarrow H$ be a game reduction. We introduce the following
definition.

\begin{definition}
$\succ _{H}$ has property $K$ if for each $i\in I$ and for each $y\in G_{i},$
there exists $z_{0}\in G_{i}$ with $z_{0}\succeq _{H}y$ such that $\{z\in
G_{i}:z\succeq _{H}z_{0}\}$ is compact.
\end{definition}

Lemma 3 generalizes the Dufwenberg-Stegeman Lemma by relaxing the continuity
assumption on the payoff functions of the game. We use the notion of
transfer upper continuity due to Tian and Zhou (1995). Note that $G$ may not
be compact.

Before stating the lemma, we define two types of discontinuous games.

\begin{definition}
The game $G$ is called
\end{definition}

i) own transfer upper continuous if $u_{i}(\cdot ,s_{-i})$ is transfer upper
continuous for each $i\in I$ and for each $s_{-i}\in G_{-i};$

ii) own transfer weakly upper continuous if $u_{i}(\cdot ,s_{-i})$ is
transfer weakly upper continuous for each $i\in I$ and for each $s_{-i}\in
G_{-i};$

\begin{lemma}
Let us assume that $G\rightarrow ^{\ast }H$ for an own-transfer weakly upper
continuous game $G$ and $\succ _{H}$ has property $K.$ If $y\succ _{H}x$ for
some $x,y\in G_{i}$ and $i\in I,$ then there exists $z^{\ast }\in H_{i}$
such that $z\nsucc _{H}z^{\ast }\succ _{H}z$ for each $z\in G_{i}.$
\end{lemma}

\textit{Proof.} Since $G\rightarrow ^{\ast }H,$ there exists a sequence of
parings $A^{t},$ $t=0,1,2...$ such that $A^{0}=G$, $A^{t}\rightarrow A^{t+1}$
$\forall t\geq 0$ and $H_{i}=\cap _{t}A_{i}^{t},$ $\forall i\in I.$

Let $Z:=\{z\in G_{i}:u_{i}(z,s_{-i})\geq u_{i}(y,s_{-i})$ $\forall s_{-i}\in
H_{-i}\}.$ According to property $K$ of $\succ _{H},$ it follows that there
exists $z_{0}\in G_{i}$ such that $z_{0}\succeq _{H}y$ and $U:=\{z\in
G_{i}:z\succeq _{H}z_{0}\}$ is compact. Since $y\succ _{H}x,$ we have that $%
H_{-i}\neq \emptyset .$ Let us define $f:U\rightarrow \mathbb{R}$ by $%
f(z)=u_{i}(z,s_{-i}^{\ast }),$ where $s_{-i}^{\ast }\in H_{-i}$ is fixed.

Since $f$ is transfer weakly upper continuous on $U,$ $f$ reaches its
maximum in $z^{\ast }\in U\subset Z.$ We note that $z^{\ast }\in Z$ and $%
y\succ _{H}x$ imply $z^{\ast }\succ _{H}x.$ If $z\succ _{H}z^{\ast }$ for
some $z\in G_{i},$ then $u_{i}(z,s_{-i})>u_{i}(z^{\ast },s_{-i})$ $\forall
s_{-i}\in H_{-i},$ implying that $z\in U$ and $f(z)>f(z^{\ast }),$
contradiction. Therefore, $z\nsucc _{H}z^{\ast }$ $\forall z\in G_{i},$ so
that $z\nsucc _{A_{t}}z^{\ast }$ $\forall z\in G_{i}$ $\forall t\geq 0$
implying that $z^{\ast }\in A_{i}^{t}$ $\forall t\geq 0.$ It follows that $%
z^{\ast }\in H_{i}.$

\begin{example}
Let $I=\{1,2\},$ $G_{1}=G_{2}=[0,2],$ $u_{i}:G_{i}\times G_{j}\rightarrow 
\mathbb{R},$
\end{example}

$u_{i}(x,y)=\left\{ 
\begin{array}{c}
1\text{ \ \ \ \ \ \ if \ \ \ \ \ }x=0; \\ 
2\text{ \ \ \ if \ \ }x\in (0,1); \\ 
x+1\text{ if }x\in \lbrack 1,2].%
\end{array}%
\right. $

Let $H=(H_{1},H_{2})$ $H_{1}=H_{2}=[0,1].$

We notice that, for each $y\in G_{2},$ $u_{i}(.,y)$ is transfer weakly upper
continuous on [0,2] and $u_{i}(.,y)$ is not upper semicontinuous at $x=0.$

We prove that $\succ _{H}$ has property $K:$

If $y=0,$ there exists $z_{0}=0$ such that $U(0)=\{z\in \lbrack
0,2]:u_{1}(z,s)\geq u_{1}(0,s)$ for each $s\in H_{2}\}=[0,2]$ is a compact
set$.$

If $y\in (0,1),$ there exists $z_{0}=\frac{3}{2}$ such that $U(\frac{3}{2}%
)=\{z\in \lbrack 0,2]:u_{1}(z,s)\geq u_{1}(\frac{3}{2},s)$ for each $s\in
H_{2}\}=[\frac{3}{2},1]$ is a compact set$.$

If $y\in \lbrack 1,2],$ there exists $z_{0}=y$ such that $U(z_{0})=\{z\in
\lbrack 0,2]:u_{1}(z,s)\geq u_{1}(z_{0},s)$ for each $s\in H_{2}\}=[y,2]$ is
a compact set$.$

We have that for any $x,y\in \lbrack 0,2]$ such that $y\succ _{H}x$, there
exists $z^{\ast }\in \lbrack 0,2]$ such that $z^{\ast }\succ _{H}x$ and $%
z\nsucc _{H}z^{\ast }\succ x$ for each $z\in H_{i}.$

If $H=G,$ we obtain the following corollary.

\begin{corollary}
Let assume that $G$ is an own-transfer weakly upper continuous game $G$ and $%
\succ _{G}$ has property $K.$ If $y\succ _{G}x$ for some $x,y\in G_{i}$ and $%
i\in I,$ then there exists $z^{\ast }\in G_{i}$ such that $z\nsucc
_{G}z^{\ast }\succ _{G}z$ for each $z\in G_{i}.$
\end{corollary}

If in the last corollary, the game $G$ is transfer upper semicontinuous and
compact, we obtain the following result.

\begin{corollary}
Let assume that $G$ is a compact, own transfer upper semicontinuous game $G.$
If $y\succ _{G}x$ for some $x,y\in G_{i}$ and $i\in I,$ then there exists $%
z^{\ast }\in G_{i}$ such that $z\nsucc _{G}z^{\ast }\succ _{G}z$ for each $%
z\in G_{i}.$
\end{corollary}

In order to obtain other generalization of Dufwenberg-Stegeman Lemma (2002),
we further define the property $\mathcal{M}$ for a function $u$.

\begin{definition}
Let $X$ be a subset of a topological space. The function $u:X\rightarrow 
\mathbb{R}\cup \{-\infty \}$ has the property $\mathcal{M}$ on $X$ if for
each $y\in X,$ $x\in $cl$\{z\in X:u(z)\geq u(y)\}\backslash \{z\in
X:u(z)\geq u(y)\},$ implies there exists $x^{\prime }\in X$ such that $%
u(x^{\prime })>u(x).$
\end{definition}

We provide an example of transfer weakly upper continuous function which
verifies the property $\mathcal{M}$.

\begin{example}
$u:[0,1]\rightarrow \mathbb{R}$, $u(x)=\left\{ 
\begin{array}{c}
1,\text{ if }x\text{ is a rational number,} \\ 
0,\text{ \ \ \ \ \ \ \ \ \ \ \ \ \ \ \ \ \ \ otherwise.}%
\end{array}%
\right. $
\end{example}

First, let $y\in Q.$ For example, let $y=\frac{1}{2},$ $u(y)=1.$ Let $%
U=\{z\in \lbrack 0,1]:u(z)\geq 1\}=[0,1]\cap \mathbb{Q}.$ The set cl$%
U=[0,1]. $ If $x\in $cl$U\backslash U=[0,1]\cap (\mathbb{R}\backslash 
\mathbb{Q}),$ $u(x)=0$ and there exists $x^{\prime }\in \lbrack 0,1]$ such
that $u(x^{\prime })=1>u(x)=0.$

The property $\mathcal{M}$ is also verified for $y\in \mathbb{R}\backslash
Q. $

Lemma 4 also generalizes Dufwenberg-Stegeman Lemma (2002).

\begin{lemma}
Let us assume that $G\rightarrow ^{\ast }H$ for a compact and own-transfer
weakly upper continuous game $G$ and for each $i\in I$ and for each $%
s_{-i}\in H_{-i},$ the function $u_{i}(\cdot ,s_{-i})$ has property $%
\mathcal{M}.$ If $y\succ _{H}x$ for some $x,y\in G_{i}$ and $i\in I,$ then
there exists $z^{\ast }\in H_{i}$ such that $z\nsucc _{H}z^{\ast }\succ
_{H}z $ for each $z\in G_{i}.\medskip $
\end{lemma}

\textit{Proof.} Since $G\rightarrow ^{\ast }H,$ there exists a sequence of
parings $A^{t},$ $t=0,1,2...$ such that $A^{0}=G$, $A^{t}\rightarrow A^{t+1}$
$\forall t\geq 0$ and $H_{i}=\cap _{t}A_{i}^{t},$ $\forall i\in I.$

Since $y\succ _{H}x,$ we have that $H_{-i}\neq \emptyset .$ For each $%
s_{-i}\in H_{-i},$ let $Z(s_{-i}):=\{z\in G_{i}:u_{i}(z,s_{-i})\geq
u_{i}(y,s_{-i})\}$ and $Z:=\cap _{s_{-i}\in H_{-i}}$cl$Z(s_{-i}).$ The set $%
Z $ is compact. Let us define $f:Z\rightarrow \mathbb{R}$ by $%
f(z)=u_{i}(z,s_{-i}^{\ast }),$ where $s_{-i}^{\ast }\in H_{-i}$ is fixed.

Since $f$ is transfer weakly upper continuous on $Z,$ $f$ attains its
maximum in $z^{\ast }\in Z.$ Each $u_{i}$ has property $\mathcal{M}$ and
then we conclude that for each $s_{-i}\in H_{-i},$ $z^{\ast }\in Z(s_{-i}).$
We have $y\succ _{H}x$ and this fact implies $z^{\ast }\succ _{H}x.$ If $%
z\succ _{H}z^{\ast }$ for some $z\in G_{i},$ then $u_{i}(z,s_{-i})>u_{i}(z^{%
\ast },s_{-i})$ $\forall s_{-i}\in H_{-i},$ implying that $z\in Z$ and $%
f(z)>f(z^{\ast }),$ contradiction. Therefore, $z\nsucc _{H}z^{\ast }$ $%
\forall z\in G_{i},$ so that $z\nsucc _{A_{t}}z^{\ast }$ $\forall z\in G_{i}$
$\forall t\geq 0$ implying $z^{\ast }\in A_{i}^{t}$ $\forall t\geq 0.$ It
follows that $z^{\ast }\in H_{i}.$

If $H=G,$ we obtain the following corollary.

\begin{corollary}
Let us assume that $G$ is a compact and own-transfer weakly upper continuous
game $G$ and for each $i\in I$ and for each $s_{-i}\in G_{-i},$ the function 
$u_{i}(\cdot ,s_{-i})$ has property $\mathcal{M}.$ If $y\succ _{G}x$ for
some $x,y\in G_{i}$ and $i\in I,$ then there exists $z^{\ast }\in G_{i}$
such that $z\nsucc _{G}z^{\ast }\succ _{G}z$ for each $z\in G_{i}.\medskip $
\end{corollary}

The next theorem is Theorem 1 in Dufwenberg and Stegeman (2002). It is the
main result concerning the existence and uniqueness of nonempty maximal
reductions of compact and continuous games.

\begin{theorem}
a) If a game $G$ is compact and own-uppersemicontinuous, then any nonempty
maximal $(\rightarrow ^{\ast })$ reduction of $G$ is the unique maximal $%
(\rightarrow ^{\ast })$ reduction of $G.$
\end{theorem}

\textit{b) If a game }$G$\textit{\ is compact and continuous, then }$G$%
\textit{\ has a unique maximal }$(\rightarrow ^{\ast })$\textit{\ reduction }%
$M;$\textit{\ furthermore, }$M$\textit{\ is nonempty, compact and continuous.%
}\medskip

We generalize the theorem above by weakening the continuity conditions on
payoff functions which describe the game model. In order to do this, we
introduce the following definition.

\begin{definition}
The reduction $G\rightarrow ^{\ast \ast }H$ is defined by the existence of
(finite or countable infinite) sequence of parings $A^{t}$ of $G,$ $%
t=0,1,2...$, such that $A^{0}=G,$ $A^{t}\rightarrow A^{t+1}$ for each $t\geq
0$ and $H_{i}=\cap _{t}A_{i}^{t}$ for each $i\in I$ and by the consistency
with the continuity of the utility functions, which means that for each $%
i\in I,$ the payoff function $u_{i}$ maintains the same continuity property
on each set $\tprod\nolimits_{i\in I}A_{i}^{t},$ $t=0,1,2...$, as it has on $%
\tprod\nolimits_{i\in I}G_{i}.$
\end{definition}

\begin{definition}
The function $u_{i}:G\rightarrow \mathbb{R}$ has the intersection property
with respect to the $i^{th}$ variable if there exists $S_{-i}\subset G_{-i}$
such that $Z_{-i}(x)=S_{-i}$ for each $x\in G_{i\text{ }}$and $Z_{i}(x)=\cap
_{s_{-i}\in S_{-i}}F_{i}(x,s_{-i}),$ where $Z(x)=\{(s_{i},s_{-i})\in
G:u_{i}(x,s_{-i})\leq u_{i}(s_{i},s_{-i})\},$ $Z_{i}(x)=$pr$_{i}Z(x),$ $%
Z_{-i}(x)=$pr$_{-i}Z(x)$ and $F_{i}(x,s_{-i})=\{s_{i}\in
G_{i}:u_{i}(x,s_{-i})\leq u_{i}(s_{i},s_{-i})\}.\medskip $
\end{definition}

\begin{example}
Let $G=[0,1]\times \lbrack 0,1]$ and let $u_{1}:G\rightarrow \mathbb{R}$ be
defined by
\end{example}

$u_{1}(x,y)=\left\{ 
\begin{array}{c}
1+x+y\text{ if }x\in Q; \\ 
x\text{ if }x\in \mathbb{R}\backslash Q.%
\end{array}%
\right. $

For each $y\in \lbrack 0,1],$ the function $u_{1}(\cdot ,y)$ is not upper
semicontinuous, but it is transfer upper continuous since, for a
neighborhood $\mathcal{N}\subset \lbrack 0,1],$ we may choose any $x^{\prime
}$ rational such that sup$\{x:x\in \mathcal{N}\}<x\prime \leq 1.$

We prove that $u_{1}$ fulfills the intersection property with respect to $x.$

We have that $Z(x)=\left\{ 
\begin{array}{c}
\{[x,1]\cap Q\}\times \lbrack 0,1]\text{ if }x\in Q; \\ 
\{[0,1]\cap Q\}\times \lbrack 0,1]\cup \{[x,1]\cap \mathbb{R}\backslash
Q\}\times \lbrack 0,1]\text{ if }x\in \mathbb{R}\backslash Q,%
\end{array}%
\right. $

$Z_{1}(x)=\left\{ 
\begin{array}{c}
\lbrack x,1]\cap Q\text{ if }x\in Q; \\ 
\{[0,x]\cap Q\}\cup \lbrack x,1]\text{ if }x\in \mathbb{R}\backslash Q%
\end{array}%
\right. $, $Z_{2}(x)=[0,1]=G_{2}$ and

$F_{1}(x,s_{2})=\left\{ 
\begin{array}{c}
\lbrack x,1]\cap Q\text{ if }x\in Q; \\ 
\{[0,x]\cap Q\}\cup \lbrack x,1]\cap \mathbb{R}\backslash Q\text{ if }x\in 
\mathbb{R}\backslash Q.%
\end{array}%
\right. $

It follows that $Z_{1}(x)=\cap _{s_{2}\in G_{2}}F_{1}(x,s_{2}).$

\begin{definition}
The game $G$ has the intersection property if $u_{i}:G\rightarrow \mathbb{R}$
has the intersection property with respect to the $i^{th}$ variable for each 
$i\in I.$
\end{definition}

\begin{theorem}
\textit{a) Let }$G$\textit{\ be an own-transfer weakly upper continuous game
which has also the intersection property, such that }$\succ _{H}$\textit{\
has property }$K$\textit{\ for every }$G\rightarrow H.$\textit{\ Then, any
nonempty maximal reduction }$G\rightarrow ^{\ast \ast }M$\textit{\ is the
unique maximal reduction.}
\end{theorem}

\textit{b) If }$G$\textit{\ is a compact, own-transfer upper continuous game
such that }$\succ _{H}$\textit{\ has property }$K$\textit{\ for every }$%
G\rightarrow H$\textit{, then it has a nonempty compact own-transfer upper
semicontinuous maximal }$(\rightarrow ^{\ast \ast })$\textit{\ reduction }$M$%
\textit{\ and this reduction is unique}$.$

\textit{Proof.} a) The proof follows the same line as Theorem 1 of
Dufwenberg and Stegeman (2002).

b) We prove that if $G\rightarrow H$ fast, then $H$ is compact and nonempty.
Since $y\succ _{G}x$ for some $x,y\in G_{i},$ then $H_{i}\neq \emptyset .$

We will show further that $H_{i}$ is compact. Let $Z(x)=\{(s_{i},s_{-i})\in
G:u_{i}(x,s_{-i})\leq u_{i}(s_{i},s_{-i})\},$ $Z_{i}(x)=$pr$_{i}Z(x)$ and $%
Z_{-i}(x)=$pr$_{-i}Z(x)=S_{-i}$ for each $x\in S_{-i}.$ Since $x\in
Z_{i}(x), $ it follows that $Z_{i}(x)\neq \emptyset .$ We first prove that $%
H_{i}=\cap _{x\in H_{i}}Z_{i}(x).$ Let us choose an arbitrary element $\ z$
of $G_{i}.$ If $z\notin Z_{i}(x),$ for each $x\in H_{i},$ then $%
u_{i}(x,s_{-i})>u_{i}(z,s_{-i})$ for each $s_{-i}\in G_{-i}.$ It follows
that $x\succ _{H}z$ and therefore, $z\notin H_{i}.$ This fact implies that $%
H_{i}\subseteq \cap _{x\in H_{i}}Z_{i}(x).$ If $z\notin H_{i},$ there exists 
$x\in G_{i}$ such that $x\succ _{G}z$ and according to Corrolary 2, it
follows that there exists $x^{\ast }\in G_{i}$ such that $x^{\ast }\succ
_{G}z$. The last assertion implies that $z\notin Z_{i}(x^{\ast })$ and
therefore, $z\notin \cap _{x\in H_{i}}Z_{i}(x).$ We have $\cap _{x\in
H_{i}}Z_{i}(x)\subseteq H_{i}$ and the equality $H_{i}=\cap _{x\in
H_{i}}Z_{i}(x)$ follows from the above assertions.

Now let us define $F_{i}(x,s_{-i})=\{s_{i}\in G_{i}:u_{i}(x,s_{-i})\leq
u_{i}(s_{i},s_{-i})\}$ and then, $Z_{i}(x)=\cap _{s_{-i}\in
S_{-i}}F_{i}(x,s_{-i}).$ Since we have the reduction $G\rightarrow ^{\ast
\ast }H,$ the function $u_{i}(.,s_{-i})$ is transfer upper continuous on $%
H_{i}$ for $s_{-i}$ fixed, and, according to Lemma 1, it follows that $\cap
_{x\in H_{i}}F_{i}(x,s_{-i})=\cap _{x\in H_{i}}$cl$F_{i}(x,s_{-i}).$

Therefore, $H_{i}=\cap _{x\in H_{i}}Z_{i}(x)=\cap _{x\in H_{i}}\cap
_{s_{-i}\in S_{-i}}F_{i}(x,s_{-i})=$

$\cap _{s_{-i}\in S_{-i}}\cap _{x\in H_{i}}F_{i}(x,s_{-i})=\cap _{s_{-i}\in
S_{-i}}\cap _{x\in H_{i}}$cl$F_{i}(x,s_{-i}),$ then $H_{i}$ is a closed set. 
$H_{i}$ is closed, $H_{i}\subset G_{i},$ $G_{i}$ is compact, then $H_{i}$ is
compact.

We consider $C(t)$ $t=0,1,...$ the unique sequence of subgames of $G$ such
that $C(0)=G$ and $C(t)\rightarrow C(t+1)$ is fast for each $t\geq 0.$ The
set $C(t)$ is compact and nonempty for each $t\geq 0.$ The game $M_{i}=\cap
_{t\geq 0}C(t)$ is compact, transfer upper semicontinuous and nonempty. We
show that $M$ is a maximal $(\rightarrow ^{\ast \ast })$-reduction of $G.$
Consider any player $i$ and $x,y\in M_{i}.$ Let $X(t):=\{s_{-i}\in
(C(t))_{-i}:u_{i}(y,s_{-i})\leq u_{i}(x,s_{-i}).$ We claim that $X(t)\neq
\emptyset .$ If not, \ for each $s_{-i}\in (C(t))_{-i}$, it follows that $%
u_{i}(y,s_{-i})>u_{i}(x,s_{-i}),$ so that $y\succ _{C(t)}x,$ contradicting $%
x\in M_{i}.$ $(C(t))_{-i}$ is compact and $\cap _{t\geq 0}(C(t))_{-i}$ is
nonempty and compact.

Let $X^{\prime }=\{s_{-i}\in \cap _{t\geq 0}(C(t))_{-i}:u_{i}(y,s_{-i})\leq
u_{i}(x,s_{-i})\}$

\ \ \ \ \ \ \ \ \ $=\{s_{-i}\in M_{-i}:u_{i}(y,s_{-i})\leq
u_{i}(x,s_{-i})\}. $

Since $M_{-i}\neq \emptyset ,$ it follows that $X^{\prime }\neq \emptyset $
and therefore $y\nsucc _{M}x$ and $M$ is maximal.

\begin{corollary}
The results also mainntain for the class of upper pseudocontinuous games.
\end{corollary}

By applying Lemma 4 , we obtain the following result.

\begin{theorem}
\textit{a) Let }$G$\textit{\ be a compact and own-transfer weakly upper
continuous game which has also the intersection property, such that for each 
}$i\in I,$\textit{\ the payoff function }$u_{i}$\textit{\ has property }$M.$%
\textit{\ Then, any nonempty maximal reduction }$G\rightarrow ^{\ast \ast }M$%
\textit{\ is the unique maximal reduction.}
\end{theorem}

\textit{b) If }$G$\textit{\ is a compact, own-transfer upper continuous game
such that for each }$i\in I,$\textit{\ the payoff function }$u_{i}$\textit{\
has property }$M$\textit{, then it has a nonempty compact own-transfer upper
semicontinuous maximal }$(\rightarrow ^{\ast \ast })$\textit{\ reduction }$%
M. $\textit{\ The reduction }$M$\textit{\ is unique.}

\section{The Mixed Strategies Case}

In Subsection 6.2 Dufwenberg and Stegeman (2002) approached the issue of
mixed strategy dominance. They distinguished between the case in which a
pure strategy is dominated by a pure strategy and the case in which it is
domintated by a mixed strategy. The main result is obtained by applying
Theorem 1 to the mixed extensions of finite games. We will extend Dufwenberg
and Stegeman's research by taking into consideration several types of
dominance relations and game reductions.

For the reader's convenience, we review here a few basic notions and
notations which deal with measurability. For an overview, please see
Parthasarathy (2005).

\subsection{Measurable spaces}

Suppose that $(G,\mathcal{G})$ is a measurable space and $H\in \mathcal{G}.$
Let us define $\mathcal{H}=\{H\cap A:A\in \mathcal{G}\}.$ $\ $Then $\mathcal{%
H}$ is a $\sigma -$algebra of subsets of $H$ and $(H,\mathcal{H})$ is a
measurable space.

\begin{definition}
Given a measurable space $(G,\mathcal{G})$ and $x\in G,$ define the
probability measure $\delta _{x}$ as
\end{definition}

$\delta _{x}(H)=\left\{ 
\begin{array}{c}
1\text{ if }x\in H; \\ 
0\text{ if \ }x\notin H%
\end{array}%
\right. $ for each $H\in \mathcal{G}.$

$\delta _{x}$ is called the Dirac measure with unit mass at $x.$

\begin{theorem}
Let $X$ be a finite set with a discrete $\sigma -$algebra. Then, every
probability $\mu $ on this measurable space can be unique represented in the
form $\mu =\tsum\nolimits_{x\in X}c_{x}\delta _{x},$ where $c_{x}\in \lbrack
0,1]$ $\forall x\in X,$ $\tsum\nolimits_{x\in X}c_{x}=1,$ thus $\mu
(E)=\tsum\nolimits_{x\in E}c_{x}$ for all $E\subset X.$
\end{theorem}

Notation If $(G,\mathcal{G})$ is a measurable space, we will denote by $%
\Delta (G)$ the set of probability measures defined on $G.$

Let $I=\{1,2,...,n\}$ and the game $G=(G_{i},u_{i})_{i\in I}.$

Assume that for each $i\in I,$ $G_{i}$ is a compact subset in a metric space 
$X$ and $u_{i}(.,s_{-i}):G_{i}\rightarrow \mathbb{R}$ is upper
semicontinuous for each $s_{-i}\in G_{i}.$

Each $u_{i}$ is measurable since it is upper semicontinuous and since it is
also bounded, it is integrable. We denote by $\Delta (G_{i})$ the set of
probability measure on the set of Borel sets on $G_{i}.$ $\Delta (G_{i})$
will be equipped with the weak topology.

\begin{theorem}
Let $G$ be a subset of a metric space. Then, $G$ is compact if and only if $%
\Delta (G)$ is compact.
\end{theorem}

A mixed strategy for player $i$ is an element $\mu _{i}\in \Delta (G_{i}).$

\begin{definition}
(Billingsley (1968), p 7). Suppose $\{\mu _{n}\}_{n\geq 1},$ $\mu _{n}$
belong to $\Delta (G),$ the set of probability measures on the Borel sets of
some compact metric space $G.$ Then "$\mu _{n}$ weakly converges to $\mu ",$
written $\mu _{n}\overset{w}{\rightarrow }\mu $ iff $\tint fd\mu
_{n}\rightarrow \tint fd\mu $ for all $f:G\rightarrow \mathbb{R},$ $f$
continuous. This topology is consistent with Prohorov metric.
\end{definition}

\begin{lemma}
(Robson 1990). Consider $u:G\rightarrow \mathbb{R}$ un upper semicontinuous
function, where $G$ is a compact metric space. It follows that $\tint ud\mu $
is upper semicontinuous in $\mu :$ $\lim \sup_{n}\tint ud\mu _{n}\leq \tint
ud\mu $ if $(\mu _{n})_{n},\mu \in \Delta (G),$ the set of probability
measures on Borel sets of $G$ and $\mu _{n}\overset{w}{\rightarrow }\mu .$

\begin{corollary}
If for the game $G,$ $u_{i}(.,s_{-i}^{\ast }):G_{i}\rightarrow \mathbb{R}$
is upper semicontinuous, then the function $V_{i}(.,s_{-i}^{\ast })$ $%
:\Delta (G_{i})\rightarrow \mathbb{R},$ defined by $V_{i}(\mu
_{i},s_{-i}^{\ast })=\tint u_{i}(\mu _{i},s_{-i}^{\ast })d\mu _{i}(s_{i})$
is upper semicontinuous.
\end{corollary}
\end{lemma}

We define the following extension of $\succ _{H}$:

\begin{definition}
Let $G\rightarrow H,$ $G=(G_{i},u_{i})_{i\in I},$ $I$ finite, $G_{i}$ is a
subset of a metric space $X$ for each $i\in I.$ Let $\Delta (G_{i})$ be the
set of probability measures on Borel sets of $G_{i}$ and $V_{i}(.,s_{-i})$ $%
:\Delta (G_{i})\rightarrow \mathbb{R},$ defined by $V_{i}(\mu
_{i},s_{-i})=\tint u_{i}(\mu _{i},s_{-i})d\mu _{i}(s_{i})$ for each $s_{-i}$
fixed. Given $x\in G_{i}$ and $\mu _{i}\in \Delta (G_{i}),$ we say that $\mu
_{i}\succ _{H}x$ if $H_{-i}\neq \emptyset $ and $V_{i}(\mu
_{i},s_{-i})>u_{i}(x,s_{-i})$ for each $s_{-i}\in H_{-i}.$
\end{definition}

\begin{lemma}
If $G\rightarrow ^{\ast }H$ for some compact and own-upper semicontinuous
game $G$ and $y\succ _{H}x$ for some $x,y\in G_{i}$ and $i\in I,$ then,
there exists $z^{\ast }\in H_{i}$ such that $\mu \nsucc _{H}z^{\ast }\succ
_{H}x$ for each $\mu \in \Delta (G_{i}).$
\end{lemma}

\textit{Proof.} The assumptions of Dufwenberg-Stegeman Lemma are fulfilled.
Then there exists $z^{\ast }\in H_{i}$ such that $z\nsucc _{H}z^{\ast }\succ
_{H}x$ for each $z\in G_{i}.$ We prove that, in addition, $\mu \nsucc
_{H}z^{\ast }$ for each $\mu \in \Delta (G_{i}).$

If $\mu \succ _{H}z^{\ast }$ for some $\mu \in \Delta (G_{i}),$ then $\tint
u_{i}(s_{i},s_{-i})d\mu _{i}(s_{i})>u_{i}(z^{\ast },s_{-i})$ for each $%
s_{-i}\in H_{-i},$ implying $\tint u_{i}(s_{i},s_{-i}^{\ast })d\mu
_{i}(s_{i})>u_{i}(z^{\ast },s_{-i}^{\ast })$ for some $s_{-i}^{\ast }$ fixed
in $H_{-i}.$.........(1)

We note that $z^{\ast }=\arg \max_{s_{i}\in Z}u_{i}(s_{i},s_{-i}^{\ast
})=\arg \max_{s_{i}\in G_{i}}u_{i}(s_{i},s_{-i}^{\ast }),$ where $Z=\{z\in
G_{i}:u_{i}(z,s_{-i})\geq u_{i}(y,s_{-i})$ $\forall s_{-i}\in H_{-i}\}.$ It
follows that $u_{i}(z^{\ast },s_{-i}^{\ast })\geq u_{i}(s_{i},s_{-i}^{\ast })
$ for each $s_{-i}\in G_{-i},$ and therefore, $u_{i}(z^{\ast },s_{-i}^{\ast
})\geq \tint u_{i}(s_{i},s_{-i}^{\ast })d\mu _{i}(s_{i}),$ relation which
contradicts (1). Therefore, $\mu \nsucc _{H}z^{\ast }$ for each $\mu \in
\Delta (G_{i}).$

\subsection{Types of dominance relations and reductions}

Let $I$ be a finite set. For each $i\in I,$ let $(G_{i},\mathcal{G}_{i})$ be
a measurable space, $H_{i}\in \mathcal{G}_{i},$ $\mathcal{H}_{i}=\{H_{i}\cap
A:A\in \mathcal{G}_{i}\}$ and $u_{i}:\tprod\nolimits_{i\in
I}G_{i}\rightarrow \mathbb{R}$ be a $\otimes \mathcal{G}_{i}-$measurable and
bounded. Let $G=\tprod\nolimits_{i\in I}G_{i}.$

\begin{definition}
We define the followings types of dominance relations.
\end{definition}

i) $\mu \succ _{\Delta (H)}m$ for $\mu ,m\in \Delta (G_{i})$ if $%
\tint\nolimits_{G_{i}\times H_{-i}}u_{i}(s_{i},s_{-i})d\mu _{1}\times
...\times d\mu _{i-1}\times d\mu \times d\mu _{i+1}\times ...\times d\mu
_{n}>\tint\nolimits_{G_{i}\times H_{-i}}u_{i}(s_{i},s_{-i})d\mu _{1}\times
...\times d\mu _{i-1}\times dm\times d\mu _{i+1}\times ...\times d\mu _{n},$ 
$\forall $ $\mu _{1}\times ...\times \mu _{i-1}\times \mu _{i+1}\times
...\times \mu _{n}\in \Delta (H_{-i}).$

ii) $\mu \succ _{H}m$ for $\mu ,m\in \Delta (G_{i})$ if $\tint%
\nolimits_{G_{i}}u_{i}(s_{i},s_{-i})d\mu
>\tint\nolimits_{G_{i}}u_{i}(s_{i},s_{-i})dm$ $\forall s_{-i}\in H_{-i}.$

iii) $\mu \succ _{H}x$ for $\mu \in \Delta (G_{i})$ and $x\in G_{i}$ if $\mu
\succ _{H}\delta _{x},$ which is equivalent with $\tint%
\nolimits_{G_{i}}u_{i}(s_{i},s_{-i})d\mu >u_{i}(x,s_{-i})$ $\forall
s_{-i}\in H_{-i}.$

We obtain the following theorem.

\begin{theorem}
With the notations above, we have the following relations amongst the former
types of dominance:
\end{theorem}

i) $\mu \succ _{\Delta (H)}m$ for $\mu ,m\in \Delta (G_{i})\Rightarrow \mu
\succ _{H}m$

To prove this fact, we take $\mu =\delta _{s_{j}}$ for $j\neq i,$ $s_{j}\in
H_{j}.$

ii) $\mu \succ _{H}m$ for $\mu ,m\in \Delta (G_{i})\Rightarrow \mu \succ
_{H}x\Leftrightarrow \mu \succ _{H}\delta _{x},$ $x\in G_{i}.$

Since $\delta _{x}\in \Delta (G_{i})$ for $x\in G_{i},$ ii) can be easily
checked.\medskip

Let us consider parings $G,H$ with the property that $H_{i}\subseteq G_{i}$
for each $i\in I.$ In addition to the game reduction used by Dufwenberg and
Stegeman (2002), we present the following ones.

\begin{definition}
i) (Gilboa, Kalai and Zemel 1990) $G\Rightarrow H$ if, for each $x\in
G_{i}\backslash H_{i},$ there exists $y\in H_{i}$ such that $y\succ _{H}x.$
\end{definition}

ii) $G\mapsto H$ if, for each $x\in G_{i}\backslash H_{i},$ there exists $%
\mu \in \Delta (G_{i})$ such that $\mu \succ _{H}x.$

iii) $G\rightrightarrows H$ if, for each $x\in G_{i}\backslash H_{i},$ there
exists $\mu \in \Delta (H_{i})$ such that $\mu \succ _{H}x.$

iv) $\Delta (G)\hookrightarrow \Delta (H)$ if, for each $m\in \Delta
(G_{i})\backslash \Delta (H_{i}),$ there exists $\mu \in \Delta (G_{i})$
such that $y\succ _{H}m.$

v) $\Delta (G)\Rrightarrow \Delta (H)$ if, for each $m\in \Delta
(G_{i})\backslash \Delta (H_{i}),$ there exists $\mu \in \Delta (H_{i})$
such that $y\succ _{H}m.\medskip $

We will need the following theorem.

\begin{theorem}
There are the following relations amongst the former types of reductions.
\end{theorem}

i) $(G\Rightarrow H)$ $\Longrightarrow $ $(G\rightarrow H)$

\ \ $(G\rightrightarrows H)$ $\Longrightarrow $ $(G\mapsto H)$

\ \ $(\Delta (G)\Rrightarrow \Delta (H))$ $\Longrightarrow $ $(\Delta
(G)\hookrightarrow \Delta (H))$

ii) $(\Delta (G)\Rightarrow \Delta (H))$ $\Longrightarrow $ $(\Delta
(G)\Rrightarrow \Delta (H))$ $\Longrightarrow $ $(G\rightrightarrows H)$

iii) $(\Delta (G)\rightarrow \Delta (H))$ $\Longrightarrow $ $(\Delta
(G)\hookrightarrow \Delta (H))$ $\Longrightarrow $ $(G\mapsto H)$

\textit{Proof.}

i) The proof is obvious.

ii) \ Suppose $(\Delta (G)\Rightarrow \Delta (H))$. It follows that $\Delta
(H_{i})\subset \Delta (G_{i})$ and for each $m\in \Delta (G_{i})\backslash
\Delta (H_{i}),$ there exists $\mu \in \Delta (H_{i})$ such that $\mu \succ
_{\Delta (H)}m.$ According to Theorem 7, it follows that $H_{i}\subset G_{i}$
for each $i\in I$ and for each $m\in \Delta (G_{i})\backslash \Delta
(H_{i}), $ there exists $\mu \in \Delta (H_{i})$ such that $\mu \succ _{H}m$ 
$\iff $ $\Delta (G)\Rrightarrow \Delta (H).$

If $m=\delta _{x}$ with $x\in G_{i}\backslash H_{i},$ we have that $%
H_{i}\subset G_{i}$ for each $i\in I$ and for each $x\in G_{i}\backslash
H_{i},$ there exists $\mu \in \Delta (H_{i})$ such that $\mu \succ _{H}x,$
which is equivalent with $G\rightrightarrows H.$

iii) The implication are true from i) and ii).

\begin{theorem}
If $G$ is a finite game, then $(\Delta (G)\Rrightarrow \Delta (H))$ $\iff $ $%
(G\rightrightarrows H).$
\end{theorem}

\textit{Proof.} The direct implication $"\Longrightarrow "$ comes from
Theorem 8, ii).

We prove $"\impliedby ".$ Let $x\in G_{i}\backslash H_{i}.$ Since $%
G\rightrightarrows H,$ there exists $\mu _{x}\in \Delta (H_{i})$ such that $%
\mu _{x}\succ _{\Delta (H)}x$ $\Longrightarrow \mu _{x}\succ _{\Delta
(H)}\delta _{x}.$

Let $m\in \Delta (G_{i})\backslash \Delta (H_{i}).$ According to Theorem 5, $%
m$ can be unique represented as a convex combination of Dirac measures $%
\delta _{x},$ $x\in G_{i}\backslash H_{i}.$ Then, there exists unique $%
c_{x}\in \lbrack 0,1],$ $\tsum\nolimits_{x\in G_{i}\backslash H_{i}}c_{x}=1$
such that $m=\tsum\nolimits_{x\in G_{i}\backslash H_{i}}c_{x}\delta _{x}.$
But, as we noted above, for each $\delta _{x}$ with $x\in G_{i}\backslash
H_{i},$ there exists $\mu _{x}\in \Delta (H_{i})$ such that $\mu _{x}\succ
_{H}\delta _{x}.$ Therefore, $\mu =\tsum\nolimits_{x\in G_{i}\backslash
H_{i}}c_{x}\mu _{x}$ is a probability measure on $G_{i}\backslash H_{i}$ and 
$\mu _{x}\succ _{H}m.$

\begin{theorem}
Let $G\rightrightarrows H.$ If there exists $x^{\ast }\in G_{i}\backslash
H_{i}$ such that $u_{i}(x^{\ast },s_{-i})\geq u_{i}(x,s_{-i})$ for each $%
x\in G_{i}\backslash H_{i}$ and $s_{-i}\in H_{-i},$ then $\Delta
(G)\Rrightarrow \Delta (H).$
\end{theorem}

\textit{Proof.} Let $x^{\ast }$ be such that $x^{\ast }\in G_{i}\backslash
H_{i}$ and $u_{i}(x^{\ast },s_{-i})\geq u_{i}(x,s_{-i})$ for each $x\in
G_{i}\backslash H_{i}$ and $s_{-i}\in H_{-i}.$ Then, $\tint\nolimits_{G_{i}%
\backslash H_{i}}u_{i}(x,s_{-i})dm\leq u_{i}(x^{\ast },s_{-i})$ for each $%
m\in \Delta (G_{i}\backslash H_{i}).$ \ \ \ \ \ \ \ \ \ \ \ \ \ \ \ \ \ \ \
\ \ \ \ \ \ \ \ \ \ \ \ \ \ \ \ \ \ (1)

Since $x^{\ast }\in G_{i}\backslash H_{i}$ and $G\rightrightarrows H,$ it
follows that there exists $\mu \in \Delta (H_{i})$ such that $\mu \succ
_{H}x^{\ast },$ that is $\tint\nolimits_{H_{i}}u_{i}(x,s_{-i})d\mu
>u_{i}(x^{\ast },s_{-i})$ for each $s_{-i}\in H_{-i}.$ \ \ \ \ \ \ \ \ \ \ \
\ \ \ \ \ \ \ \ \ \ \ \ (2)

From 1) and 2), it follows that for $m\in \Delta (G_{i})\backslash \Delta
(H_{i}),$ there exists $\mu \in \Delta (H_{i})$ such that $%
\tint\nolimits_{H_{i}}u_{i}(x,s_{-i})d\mu >\tint\nolimits_{G_{i}\backslash
H_{i}}u_{i}(x,s_{-i})dm$ for each $s_{-i}\in H_{-i},$ that is $\mu \succ
_{H}m.$ Therefore, $\mu \succ _{H}m.$

\begin{corollary}
Let $G\mapsto H.$ If there exists $x^{\ast }\in G_{i}\backslash H_{i}$ such
that $u_{i}(x^{\ast },s_{-i})\geq u_{i}(x,s_{-i})$ for each $x\in
G_{i}\backslash H_{i}$ and $s_{-i}\in H_{-i},$ then $\Delta
(G)\hookrightarrow \Delta (H).$
\end{corollary}

\subsection{Dufwenberg-Stegeman-like Lemma}

We study first the case of the game reduction $G\mapsto H$.

\begin{lemma}
In the case of a finite game, Lemma Dufwenberg-Stegeman remains true for the
game reduction $G\mapsto H.$
\end{lemma}

\textit{Proof.} If $G\mapsto ^{\ast }H,$ then $\Delta (G)\hookrightarrow
^{\ast }\Delta (H),$ according to Theorem 9 and Theorem 8. Let $\mu ^{\prime
}\succ _{H}x$ for some $x\in G_{i}$ and $\mu ^{\prime }\in \Delta (G_{i}).$
Then, $\mu ^{\prime }\succ _{H}\delta _{x}.$ By applying Lemma
Dufwenberg-Stegeman to $\Delta (G)\hookrightarrow ^{\ast }\Delta (H),$ we
obtain that there exists $\mu ^{\ast }\in \Delta (H_{i})$ such that $\mu
\ngtr _{H}\mu ^{\ast }\succ _{H}\delta _{x}$ for each $\mu \in \Delta
(G_{i}).$ Therefore, there exists $\mu ^{\ast }\in \Delta (H_{i})$ such that 
$\mu \nsucc _{H}\mu ^{\ast }\succ _{H}x$ for each $\mu \in \Delta (G_{i}).$

We also obtain the next result concerning the game reduction $G\mapsto
^{\ast }H.$

\begin{lemma}
Let $I$ be a finite set. For each $i\in I,$ let $G_{i}$ be a compact subset
of a metric space $X$ considered with its borelian sets, $H_{i}\subset G_{i}$
and $u_{i}:\tprod\nolimits_{i\in I}G_{i}\rightarrow \mathbb{R}_{+}$
uppersemicontinuous in each argument. Let $G\mapsto ^{\ast }H$ and suppose
that $\Delta (G)\hookrightarrow ^{\ast }\Delta (H).$ If $\mu ^{\prime }\succ
_{H}x$ for some $x\in G_{i}$ and $\mu ^{\prime }\in \Delta (G_{i}),$ $i\in
I, $ then there exists $\mu ^{\ast }\in \Delta (H_{i})$ such that $\mu
\nsucc _{H}\mu ^{\ast }\succ _{H}x$ for each $\mu \in \Delta (G_{i}).$
\end{lemma}

\textit{Proof.} According to Theorem 6, if $G$ is compact, $\Delta (G)$ is
also compact. According to Corollary 2, if $u_{i}(.,s_{-i})$ is upper
semicontinuous for each $s_{-i}\in G_{-i},$ then $V_{i}(.,s_{-i})$ is also
upper semicontinuous on $\Delta (G_{i})$ for each $s_{-i}\in G_{-i},$ where $%
V_{i}(\mu ,s_{-i})=\tint\nolimits_{G_{i}}u_{i}(s,s_{-i})d\mu .$

Since $G\mapsto ^{\ast }H,$ there exists a sequence of parings $A^{t},$ $%
t=0,1,2...$ such that $A^{0}=G$, $A^{t}\mapsto A^{t+1}$ $\forall t\geq 0$
and $H_{i}=\cap _{t}A_{i}^{t},$ $\forall i\in I.$ Let $Z=\{\mu \in \Delta
(G_{i}):V_{i}(\mu ,s_{-i})\geq V_{i}(\mu ^{\prime },s_{-i})$ for each $%
s_{-i}\in H_{-i}\}.$ $Z$ is a nonempty set, since $\mu ^{\prime }\in Z.$ $Z$
is also closed (as intersection of the closed sets $Z(s_{-i})=\{\mu \in
\Delta (G_{i}):V_{i}(\mu ,s_{-i})\geq V_{i}(\mu ^{\prime },s_{-i}),$ $%
s_{-i}\in H_{-i}\}$) and therefore compact. Let us define $f:Z\rightarrow 
\mathbb{R},$ $f(\mu )=V_{i}(\mu ,s_{-i}^{\ast })$ for $s_{-i}^{\ast }\in
H_{-i}$ fixed. The function $f$ is uppersemicontinuous and it reaches its
maximum on the compact set $Z.$ Denote by $\mu ^{\ast }=\arg \max_{\mu \in
Z}f(\mu ).$

It follows that, there exists $\mu ^{\ast }\in \Delta (G_{i})$ such that $%
\mu \ngtr _{H}\mu ^{\ast }\succ _{H}x$ for each $\mu \in \Delta (G_{i}).$
Therefore, $\mu \ngtr _{H}\mu ^{\ast }$ for each $\mu \in \Delta (G_{i})$
and then $\mu \ngtr _{A^{t}}\mu ^{\ast }$ for each $\mu \in \Delta (G_{i})$
and $t\geq 0$. Since $\Delta (G)\hookrightarrow ^{\ast }\Delta (H),$ we
conclude that $\mu ^{\ast }\in \Delta (A^{t})$ for each $t\geq 0,$ and
therefore, $\mu ^{\ast }\in \Delta (H_{i}).\medskip $

For $H=G,$ we obtain the following corollary.

\begin{corollary}
Let $I$ be a finite set. For each $i\in I,$ let $G_{i}$ be a compact subset
of a metric space $X$ considered with its borelian sets and $%
u_{i}:\tprod\nolimits_{i\in I}G_{i}\rightarrow \mathbb{R}_{+}$ upper
semicontinuous in each argument. If $\mu ^{\prime }\succ _{G}x$ for some $%
x\in G_{i}$ and $\mu ^{\prime }\in \Delta (G_{i}),$ $i\in I,$ then there
exists $\mu ^{\ast }\in \Delta (G_{i})$ such that $\mu \nsucc _{G}\mu ^{\ast
}\succ _{G}x$ for each $\mu \in \Delta (G_{i}).$
\end{corollary}

\begin{corollary}
Lemma 8 is true if, instead of having the assumption $\Delta
(G)\hookrightarrow ^{\ast }\Delta (H),$ we have the following one: there
exists $x^{\ast }\in G_{i}\backslash H_{i}$ such that $u_{i}(x^{\ast
},s_{-i})\geq u_{i}(x,s_{-i})$ for each $x\in G_{i}\backslash H_{i}$ and $%
s_{-i}\in H_{-i}.$
\end{corollary}

\textit{Proof.} The proof of the corollary comes from Theorem 6.

\subsection{Existence and uniqueness of maximal reductions}

The main result of Section 4 is Theorem 11.

\begin{theorem}
Let $G=(I,(G_{i})_{i\in I},(u_{i})_{i\in I})$ be a strategic game such that $%
I$ is a finite set and for each $i\in I,$ $G_{i}$ is a nonempty compact
subset of a metric space, $u_{i}:\tprod\nolimits_{i\in I}G_{i}\rightarrow 
\mathbb{R}$ is upper semicontinuous in each argument and for each $G\mapsto
H,$ $\Delta (G)\hookrightarrow \Delta (H)$ (or there exists $x^{\ast }\in
G_{i}\backslash H_{i}$ such that $u_{i}(x^{\ast },s_{-i})\geq u_{i}(x,s_{-i})
$ for each $x\in G_{i}\backslash H_{i}$ and $s_{-i}\in H_{-i}).$ Then, $G$
has a unique nonempty maximal $(\mapsto ^{\ast })$ reduction $M$ and $M$ is
nonempty, compact and upper semicontinuous.\medskip 
\end{theorem}

\textit{Proof.} The game $(I,(\Delta (G_{i}))_{i\in I},(V_{i})_{i\in I})_{%
\text{ }}$is also compact and own-upper-semicontinuous. According to Lemma
8, we have that if $\mu ^{\prime }\succ _{H}x$ for some $x\in G_{i}$ and $%
\mu ^{\prime }\in \Delta (G_{i}),$ $i\in I,$ then there exists $\mu ^{\ast
}\in \Delta (H_{i})$ such that $\mu \nsucc _{H}\mu ^{\ast }\succ _{H}x$ for
each $\mu \in \Delta (G_{i}).$ The set $\Delta (H_{i})$ is nonempty since $%
H_{i}$ is nonempty.

The proof of the uniqueness of $M$ follows the same line as the proof of
Theorem 1a) of Dufwenberg-Stegeman.

Now we are proving that, if $G$ is compact and own-upper semicontinuous and $%
G\mapsto H$ fast, then $H$ is compact and nonempty.

Choose $i\in I$ such that $H_{i}\neq G_{i}.$ Since $\mu >_{G}x$ \ for some $%
x\in G_{i},$ $\mu \in \Delta (G_{i}),$ according to Corollary 4, we have
that $H_{i}\neq \emptyset .$ It remains to show that $H_{i}$ is compact.
Choose $\mu \in \Delta (H_{i})$ and let Let $Z(\mu )=\{(s_{i},s_{-i})\in
G:V_{i}(\mu ,s_{-i})\leq u_{i}(s_{i},s_{-i})\},$ $Z_{i}(\mu )=$pr$_{i}Z(\mu
) $ and $Z_{-i}(\mu )=$pr$_{-i}Z(\mu ).$

The set $Z_{i}(\mu )$ is nonempty. In order to prove this fact, we will
assume the opposite: $Z_{i}(\mu )=\emptyset .$ In this case, $V_{i}(\mu
,s_{-i})>u_{i}(s,s_{-i})$ for each $s\in G_{i}$ and for each $s_{-i}\in
G_{-i},$ and it follows that $\mu >_{G}s$ for each $s\in G_{i}.$ We can
conclude that $\mu \succ _{H}s$ for each $s\in G_{i},$ and, since $%
G\rightrightarrows H$ fast, we have that $($for each $s\in G_{i}\Rightarrow
s\notin H_{i})$ and, then, $H_{i}$ is an empty set, and we reached a
contradiction.

Now let us define $F_{i}(\mu ,s_{-i})=\{s_{i}\in G_{i}:V_{i}(\mu
,s_{-i})\leq u_{i}(s_{i},s_{-i})\}$ and then, $Z_{i}(\mu )=\cap _{s_{-i}\in
Z_{-i}(\mu )}F_{i}(x,s_{-i}).$

Since $u_{i}(.,s_{-i})$ is upper semicontinuous for each $s_{-i}\in G_{-i},$
we have that $Z_{i}(\mu )$ is closed as being an intersection of closed
subsets. We will show that $H_{i}=\cap _{\mu \in \Delta (H_{i})}Z_{i}(\mu ).$

Let us consider $x\in G_{i}.$ For any $\mu \in \Delta (H_{i}),$ if $x\notin
Z(\mu ),$ we have that $V_{i}(\mu ,s_{-i})>u_{i}(x,s_{-i})$ for each $%
s_{-i}\in G_{-i}$ and, therefore $\mu >_{G}x.$ Then $x\notin H_{i}$ and $%
H_{i}\subseteq \cap _{\mu \in \Delta (H_{i})}Z(\mu ).$

If $x\notin H_{i},$ then $\mu >_{G}x$ for some $\mu \in \Delta (G_{i})$ and
Lemma 8 implies that there exists $\mu ^{\ast }\in \Delta (G_{i})$ such that 
$\mu ^{\ast }>_{G}x$ and therefore, $x\notin Z(\mu ^{\ast })$ and we can
conclude that $x\notin \cap _{\mu \in \Delta (H_{i})}Z(\mu ).$ Therefore, $%
H_{i}\supseteq \cap _{\mu \in \Delta (H_{i})}Z(\mu ).$

The equality $H_{i}=\cap _{\mu \in \Delta (H_{i})}Z_{i}(\mu )$ holds and,
since $Z_{i}(\mu )$ is closed for all $\mu ,$ $H_{i}$ is also closed and
therefore compact.

Let $C\left( t\right) ,$ $t=0,1,...$ denote the unique sequence of subgames
of $G$ such that $C(0)=G$ and $C(t)\mapsto C(t+1)$ is fast for each $t\geq
0. $ We have that $C(t)$ is compact and nonempty for each $t\geq 0.$ It
follows that $M_{i}=\cap _{t\in C(t)_{i}}$ is compact, nonempty for each $%
i\in I.$

We will show that $M$ is a maximal $(\mapsto ^{\ast })$ reduction of $G.$
Let $x\in M_{i},$ $\mu \in \Delta (M_{i}).$ Let $X(t)=\{s_{-i}\in
C(t)_{-i}:V_{i}(\mu ,s_{-i})\leq u_{i}(x,s_{-i})\}.$ If $X(t)=\emptyset $
for each $t$ such that $C(t)\neq M,$ then $\mu \succ _{C(t)}x,$
contradiction. Therefore, $X(t)\neq \emptyset .$ The set $C(t)_{-i}$ is
compact for each $t$ such that $C(t)\neq M.$ Then, $M_{-i}\neq \emptyset $
and it follows that the set $X=\{s_{-i}\in M_{-i}:V_{i}(\mu ,s_{-i})\leq
u_{i}(x,s_{-i})\}$ is nonempty. We conclude that $\mu \nsucc _{M}x.\medskip $

\section{Concluding remarks}

We identified a class of discontinuous games for which the iterated
elimination of strictly dominated strategies produce a unique maximal
reduction that is nonempty. We also provided conditions under which order
independence remains valid for the case that the pure strategies are
dominated by mixed strategies. Our results expel M. Dufwenberg and M.
Stegeman's idea in [6] that 'the proper definition and the role of iterated
strict dominance is unclear for games that are not compact and continuous'.
G. Tian and J. Zhou's notion of transfer upper continuity proved to be a
suitable assumption for the payoff functions of a game in order to obtain
our results. Their Weierstrass-like theorem for transfer weakly upper
continuous functions defined on a compact set was the key of the proofs of
Lemma 3 and Lemma 4. We can conclude and emphasize that, even outside the
continuous class of games, the iterated elimination of strictly dominated
strategies remains an interesting procedure.\medskip 

We thank to Professor Krzysztof Apt for the precious ideas under which this
paper has been developed and for the hospitality he proved during our
postdoctoral stage at the Institute for Logic, Language and Computation from
the University of Amsterdam in the summer of 2012.

\end{document}